\documentclass[12pt]{article}

\usepackage{hyperref}

\usepackage{amsmath}
\usepackage{amssymb}
\usepackage{amsfonts}
\usepackage{amsopn}
\usepackage[all]{xy}
\usepackage{amsthm}
\usepackage{amscd}

\swapnumbers
\theoremstyle{plain}
\newtheorem{thm}{Theorem}[section]
\newtheorem{lem}[thm]{Lemma}

\newtheorem{prop}[thm]{Proposition}

\theoremstyle{definition}
\newtheorem{ntt}[thm]{}

\newcommand{\af}{\mathbb{A}}   
\newcommand{\gm}{{\mathbb{G}_m}} 

\newcommand{\m}{\mathfrak{m}}             
\newcommand{\spr}[1]{\langle #1 \rangle}  

\DeclareMathOperator{\spec}{Spec}   

\title{On the norm principle for quadratic forms} 

\author{M.Ojanguren, I.Panin, K.Zainoulline}

\begin{document}

\maketitle

\begin{abstract}
We prove a version of Knebusch's Norm Principle
for finite \'etale extensions
of semi-local Noetherian domains with infinite residue fields
of characteristic different from 2. As an application
we prove  Grothendieck's conjecture
on principal homogeneous spaces for the
spinor group of a quadratic space.

\


\end{abstract}

\section{Introduction}

Let $E/F$ be a finite field extension
and $q:V \to F$  a  quadratic space over $F$.
Let $D_q(F)\subset F^*$ (resp. $D_q(E)$) be the subgroup generated
by the set of non-zero elements of the field $F$ (resp. $E$)
represented by the form $q$ (resp. $q_E=q\otimes_F E$).
The well-known Knebusch's Norm Principle 
for quadratic forms over fields
\cite{KN}, \cite[VII.5.1]{La}
says that $N^E_F(D_q(E))\subset D_q(F)$,
where $N^E_F: E^*\to F^*$ is the norm map.
The goal of the present article is to show that the
Knebusch's Norm Principle holds 
for finite \'etale extensions of semi-local Noetherian domains 
with infinite residue fields of characteristic different from 2 
(see Theorem~\ref{normprrel}).
Previously, the Norm Principle for quadratic spaces
over semi-local rings was proved 
for $F$ of characteristic $0$ in \cite{Za}.
As an application we prove  Grothendieck's conjecture
on principal homogeneous spaces 
for the spinor group of a quadratic space 
(see Theorem~\ref{groth}).

This article is organized as follows.
Section~2 is devoted to some preliminary results.
In Section~3 we prove the Norm Principle 
for quadratic spaces over local rings.
Finally, we prove Grothendieck's conjecture (Sect.~\ref{grothconj}).

\subparagraph{Acknowledgments}
The second author thanks very much for the support the
RTN-Network HPRN-CT-2002-00287, 
the grant of the years 2001 to 2003 
of the ``Support Fund of National Science" 
at the Russian Academy of Science, 
the RFFI-grant 03-01-00633a and  
the Swiss National Science Foundation. 
The last author is also grateful to SFB478, M\"unster and 
AvH Foundation for hospitality and support.

\section{Preliminary Lemmas}\label{casefields}

In the present section we state and prove 
few auxiliary results which are the main tools
in the proof of Theorem \ref{normprrel}.

\begin{ntt}
Let $F$ be an infinite field 
of characteristic different from 2.
Let $E$ be a finite \'etale $F$-algebra of degree $n$.
Let $(V,q)$ be a quadratic space over $F$ of rank $m$.
Let $(V_E,q_E)$ be the base change of $(V,q)$
via the extension $E/F$, i.e.,
$V_E=V\otimes_F E$ and $q_E=q\otimes_F E$.
Sometimes we identify the vector space $V_E$
with the set $\af^m_E(E)$ of $E$-points
of the affine space $\af^m_E$.
\end{ntt}

\begin{ntt}\label{notetsep}
Since $E$ is \'etale over $F$, 
there exists an element $\alpha \in E^*$
such that
the powers $1,\alpha,\ldots,\alpha^{n-1}$
form a basis of the $F$-vector space $E$.
Such an element is called {\it primitive}.
In other words, $E$ can be written as 
the quotient $F[t]/(f_\alpha(t))$ 
of the polynomial ring $F[t]$ modulo 
the ideal generated by a monic separable polynomial
$f_\alpha(t)$ of degree $n$.
In this case $\alpha$ is identified with 
the image of $t$ by means of the quotient map
$F[t]\to F[t]/(f_\alpha(t))$.
The polynomial $f_\alpha(t)$ is called 
a {\it minimal polynomial} of $\alpha$.
\end{ntt}

\begin{ntt}\label{notet}
Observe that the subset of primitive elements of $E$
is big enough in the following sense:
Let $\beta=b_0 + b_1\alpha +\cdots+b_{n-1}\alpha^{n-1}$ with 
$b_0,\dots,b_{n-1}\in F$ be an element of $E$. 
Then, for every integer $i$ with $0\leq i\leq n-1$, 
$$
\beta^i=b_0^{(i)}+b_1^{(i)}\alpha+\cdots+b_{n-1}^{(i)}\alpha^{n-1},
\;\text{with}\; b_j^{(i)}\in F.
$$
Clearly, each $b_j^{(i)}$ is a polynomial in $b_0,\dots,b_{n-1}$. 
The condition that $\beta$ be primitive is 
the non-vanishing of the determinant 
of $\big(b_j^{(i)}\big)$ and can, thus, 
be expressed by the non-vanishing of a polynomial 
$d(b_0,\dots,b_{n-1})$. 
We define an open subset $P$ of the Weil restriction
$R_{E/F}(\gm_{,E})$ of the multiplicative group $\gm_{,E}$ by
$$
P=\{(b_0,\dots,b_{n-1})\;\vert\; d(b_0,\dots,b_{n-1})\neq 0\}\;.
$$
By  definition, the set $P(F)$ of the $F$-points of $P$
is the subset of primitive invertible elements of $E$.
\end{ntt}

\begin{lem}\label{ratlem}
Let $X$ be an $F$-rational variety.
Then any non-empty open subset $U$ of $X$
has a rational point.
\end{lem}

\begin{proof} The proof follows from the fact that
$F$ is infinite.
\end{proof}

 From this point onwards we will use the
following notations and terminology.
 
\begin{ntt}
{\it Assume there is given an element $\alpha\in P(E)$ and
a vector $v\in V_E$ such that $\alpha q(v)=1$.
Denote by $f_\alpha(t)\in F[t]$
the minimal polynomial of $\alpha$.}
\end{ntt}

\begin{ntt}
Let $Q$ be the quadric over $E$ given by
$$
Q=\{\omega\in \af^m_E \mid \alpha q(\omega)=1\}.
$$
Since the element $v$ is an $E$-point of $Q$,
the variety $Q$ is $E$-rational.
Set $Y=R_{E/F}(Q)$ to be the Weil restriction of $Q$ (see \cite{DG}).
The variety $Y$ is a closed subvariety of dimension $(m-1)n$
of the affine space $R_{E/F}(\af^m_E)=\af^{nm}_F$
and we have the bijection $Y(F)=Q(E)$ 
between the sets of $F$-points of $Y$ and $E$-points of $Q$.
Since $Q$ is rational over $E$, $Y$ is rational over $F$.
\end{ntt}

\begin{ntt}\label{phidef}
In order to define the next variety $U$
we identify $R_{E/F}(\af^m_E)$ with the
affine space $M_{n,m}$ of $n\times m$-matrices over $F$
by choosing  $\{1,\alpha,\ldots,\alpha^{n-1}\}$
as a basis of the vector space $E$ over $F$.
Thus, a vector
$\omega=(\omega_{0},\ldots,\omega_{m-1})\in
R_{E/F}(\af^m_E)(F)$
corresponds to the matrix
$(\omega_{i,j})_{i,j=0}^{n-1,m-1}$,
where the entries $\omega_{i,j}\in F$ are defined by
$\omega_{j}=\sum_{i=0}^{n-1}\omega_{i,j}\alpha^{i}$.
For any 
$\omega=(\omega_=,\dots,\omega_{n-1})$ we define 
$\omega(t)=(\omega_0(t),\dots,\omega_{m-1}(t))$, 
where 
$\omega_j(t)=\sum_{i=0}^{n-1}\omega_{i,j}t^i$.
Note that each 
$\omega_j(t)$ is of degree at most $n-1$.
We define $U$ to be the open subset
of the affine space $R_{E/F}(\af^m_E)$ defined by
$$
U=\{\omega\in M_{n,m} \mid
\Phi_{\omega}(t)=tq(\omega(t))-1
\;\text{is separable of degree}\; 2n-1\}.
$$
Clearly, the coefficients of the
polynomial $\Phi_{\omega}(t)\in F[t]$
depend on the choice
of the isomorphism $R_{E/F}(\af^m_E)\cong M_{n,m}$,
i.e., they depend on the choice of $\alpha$.
\end{ntt}

\begin{lem}\label{myclm}
The open subset $U\cap Y$ of $Y$ contains an $F$-point.
\end{lem}

\begin{proof}
Since $Y$ is $F$-rational  and
$U \cap Y$ is open in $Y$,
by Lemma \ref{ratlem} it is enough to show that
$U\cap Y$ contains a point $\rho$ 
over the algebraic closure
$\bar{F}$ of the field $F$.
Let $g(t)$ be a polynomial over $F$ such that
\begin{enumerate}
\item $\deg g(t)=n-1$
\item $g(t)$ is coprime with the polynomial
$f_\alpha(t)$
\item $-g(0)f_\alpha(0)=1$
\item $g(t)$ is separable
\end{enumerate}

Choose a hyperbolic plane $\mathbb{H}$ in
the quadratic space $(V_{\bar{F}},q_{\bar{F}})$
over $\bar{F}$. 
Choose a basis $\{e_1,e_2\}$ of $\mathbb{H}$ 
such that $q(e_1)=q(e_2)=0$ and $q(e_1+e_2)=1$.
Let $g_1(t)$, $g_2(t)$
be polynomials over $\bar{F}$
of degree $n-1$ such that
$$
tg_1(t)g_2(t)=g(t)f_\alpha(t)+1.
$$
Set $\rho(t)=g_1(t)e_1+g_2(t)e_2$ to be a vector of
polynomials over $\bar{F}$ which can be identified
with an $\bar{F}$-point of
$R_{E/F}(\af^m_E)$ following  \ref{phidef}.
Then we have
$$
tq(\rho(t))=t(q(e_1)g_1(t)^2+q(e_1+e_2)g_1(t)g_2(t)+q(e_2)g_2(t)^2)=
$$
$$
=tg_1(t)g_2(t)=g(t)f_\alpha(t)+1
$$
and
\begin{itemize}
\item $\alpha q(\rho(\alpha))=g(\alpha)f_\alpha(\alpha)+1=1$
\item $\Phi_{\rho}(t)=tq(\rho(t))-1=g(t)f_\alpha(t)$
is separable of degree $2n-1$.
\end{itemize}
Hence $\rho \in (U\cap Y)(\bar{F})$ is the desired point.
\end{proof}

\begin{ntt}\label{setw}
Define one more open subset of $Y$.
For that consider a closed subset
$Z\subset \af^m_E$ defined by
$$
Z=\{\omega\in \af^m_E\mid
\spr{v,\omega} - 1\in \{0\}_E\},
$$
where $\spr{\,,\,}: V_E\times V_E\to E$ is the bilinear form
associated with the quadratic form
$\alpha q_E$ and $\{0\}_E$
the image of the zero section $\spec E\to \af^m_E$ of $\af^m_E$.
Set $W=R_{E/F}(Q\smallsetminus Z)$.
Passing to the algebraic closure $\bar{F}$ we see that
$W$ is a non-empty open subset of $Y=R_{E/F}(Q)$.
The set $W(F)$ of $F$-points of $W$
consists of all $w\in Q(E)$
satisfying the condition
$\spr{v,\omega} - 1 \in E^*$.
\end{ntt}

\begin{lem}\label{uslem}
There exists $\omega'\in V_E$ such that
\begin{enumerate}
\item[\rm (i)] $\alpha q(\omega')=1$
\item[\rm (ii)] $\spr{v,\omega'} - 1\in E^*$
\item[\rm (iii)] the polynomial $\Phi_{\omega'}(t)$ is separable of 
degree $2n-1$.
\end{enumerate}
\end{lem}

\begin{proof}
By Lemma \ref{myclm} the set $U\cap Y$ 
is non-empty open in $Y$.
The set $W$ defined in \ref{setw}
is also non-empty open in $Y$.
Since the variety $Y$ is irreducible, the set
$W\cap U\cap Y$ is non-empty open in $Y$.
Since $Y$ is $F$-rational, the set $W\cap U\cap Y$
contains an $F$-point.
Recall that $Y(F)=Q(E)\subset V_E$.
Let $\omega'\in (W\cap U\cap Y)(F)\subset V_E$.
We claim that $\omega'$ satisfies (i) to (iii).
In fact,  property (i) holds because
$\omega'\in Q(E)$,  property (ii) holds because
$\omega'\in W(F)=(Q\smallsetminus Z)(E)$ and
property (iii) holds because $\omega'\in U(F)$.
\end{proof}

\section{The Norm Principle}\label{caserings}

\begin{ntt}
Let $R$ be a semi-local Noetherian domain
with infinite residue fields
of characteristic different from 2
(in this case $\frac{1}{2}\in R$).
Let $S/R$ be a finite \'etale $R$-algebra
(not necessarily a domain).
Let $(V,q)$ be a quadratic space  of rank $m$ over $R$.
Let $(V_S,q_S)$ be the base change of $(V,q)$
via the extension $S/R$.
Let  $D_q(R)$ (resp. $D_q(S)$)
be the group generated by the invertible
elements of $R$ (resp. $S$) represented by the form $q$.
\end{ntt}

The goal of the present section is to prove
the following

\begin{thm} \label{normprrel}
There is an inclusion
of the subgroups of $R^*$
$$
N^S_R(D_q(S)) \subset D_q(R),
$$
where $N^S_R: S^*\to R^*$ is the norm map
for the finite \'etale extension $S/R$.
\end{thm}

\begin{ntt}
For simplicity we will consider only the case
of local $R$.
By a variety over $R$ we will mean
a reduced separated scheme
of finite type over $\spec R$.
From this point onwards, by ``bar''
we mean the reduction modulo
the maximal ideal $\m$ of $R$.
So that $\overline{S}=S/\m S$.
We will write $F$ for $\overline{R}$ and
$E$ for $\overline{S}$.
So $E/F$ is a finite \'etale algebra.
\end{ntt}

To prove Theorem \ref{normprrel} we need
the following auxiliary results.

\begin{lem}\label{imlem}
Let $(S^m,\phi)$ be a quadratic space over $S$ and
let $\spr{\,,\,}: S^m\times S^m \to S$
be the associated bilinear form.
Let $v\in S^m$ be such that $\phi(v)=1$.
Let $\omega'\in E^m$ be such that
$\overline{\phi}(\omega')=\bar{1}$ and
$\spr{\overline{v},\omega'} - \bar{1} \in E^*$ (a unit).
Then there exists $\omega\in S^m$
satisfying the conditions
\begin{enumerate}
\item[\rm (i)] $\phi(\omega)=1$
\item[\rm (ii)] $\overline{\omega}=\omega'$ in $E^m$.
\end{enumerate}
\end{lem}

\begin{proof}
If $\tilde\omega$ is a lift of $\omega'$ we have 
$\phi(\tilde\omega)=1+h$ and 
$\spr{ v,\tilde\omega} - 1 = u$ with 
$h\in \m$ and $u\in R^*$. 
Putting
$$
\omega={\lambda v+\tilde\omega 
\over \lambda +1}
$$
we find 
$$
\phi(\omega)={\lambda^2+2\lambda(1+u) 
+1+h\over (\lambda+1)^2}\;.
$$
Thus, for $\lambda=-h/2u$ we have 
$\phi(\omega)=1$ and since $h\in \m$, 
$\omega$ is a lift of $\omega'$.
\end{proof}

\begin{ntt}\label{seppol}
Recall that a polynomial $f(t)\in R[t]$
is said to be separable if the quotient ring $R[t]/(f(t))$
is a finite \'etale extension of $R$.
Since $R$ is local,
a polynomial $f$ is separable iff its reduction
$\overline{f}$ modulo the maximal ideal $\m$
is a separable polynomial over $\overline{R}$.
\end{ntt}

\begin{ntt}
Similar to \ref{notetsep},
since, $S$ being \'etale over $R$, there exists
an element $\alpha \in S^*$
such that
the powers $1,\alpha,\ldots,\alpha^{n-1}$
form a basis of the free $R$-module $S$.
As in the case of fields
such an element is called {\it primitive}.
In other words, $S$ can be written as a
quotient $R[t]/(f_\alpha(t))$ of the polynomial ring
$R[t]$, where $f_\alpha(t)$ is a monic  separable polynomial
called the {\it minimal polynomial} for $\alpha$.
In this case $\alpha$ is identified with the image
of $t$ by means of the quotient map
$R[t]\to R[t]/(f_\alpha(t))$.
\end{ntt}

\begin{ntt}\label{notetrel}
As in \ref{notet} consider a primitive element
$\alpha$ of the extension $S/R$.
An element $\beta\in S^*$ can be written as $\beta=b_0 + b_1\alpha 
+\cdots+b_{n-1}\alpha^{n-1}$ in $S^*$ with $b_0,\dots,b_{n-1}\in R$. 
Then, for every integer $i$ with $0\leq i\leq n-1$, 
$$
\beta^i=b_0^{(i)}+b_1^{(i)}\alpha+\cdots+b_{n-1}^{(i)}\alpha^{n-1}
\; \text{with}\; b_j^{(i)}\in R.
$$
Clearly, each $b_j^{(i)}$ is a polynomial in 
$b_0,\dots,b_{n-1}$. The condition that $\beta$ be primitive is the 
non-vanishing of the determinant of $\big(b_j^{(i)}\big)$ and can thus 
be expressed by the non-vanishing of a polynomial 
$d(b_0,\dots,b_{n-1})$. We define an open subset $P$ of 
$R_{S/R}(\gm_{,S})$ by
$$
P=\{(b_0,\dots,b_{n-1})\;\vert\; d(b_0,\dots,b_{n-1})\in R^*\}\;.
$$
Observe that $P(R)$ is the set
of primitive invertible elements
of $S$.
\end{ntt}

\begin{ntt}\label{phidefrel}
Following \ref{phidef}, for a given primitive element
$\alpha$ and any vector $\omega\in V_S$
we define the polynomial
$\Phi_{\omega}(t)\in R[t]$ as follows.
We identify the free $S$-module $V_S$
with the set of matrices $M_{n,m}(R)$
using  $\{1,\alpha,\ldots,\alpha^{n-1}\}$
as a basis of $S$ over $R$.
Thus, a vector
$\omega=(\omega_{0},\ldots,\omega_{m-1})\in V_S$
corresponds to the matrix
$(\omega_{i,j})_{i,j=0}^{n-1,m-1}$,
where the entries $\omega_{i,j}\in R$ are defined by
$\omega_{j}=\sum_{i=0}^{n-1}\omega_{i,j}\alpha^{i}$.
We set
$$
\Phi_{\omega}(t)=tq(\omega(t))-1 \in R[t],
$$
where $\omega(t)=(\omega_0(t),\dots,\omega_{m-1}(t))$
is the vector of $m$ polynomials
$\omega_j(t)=\sum_{i=0}^{n-1}\omega_{i,j}t^i$
of degree at most $n-1$.
Clearly, the coefficients of the
polynomial $\Phi_{\omega}(t)$ depend on the choice
of $\alpha$.
\end{ntt}

\begin{lem}\label{keylem}
Let $\alpha\in S$ be a primitive invertible
element of the finite \'etale extension
$S/R$. Consider the quadratic space
$(V_S,\alpha q_S)$ over $S$.
Assume $\alpha q(v)=1$ for a vector $v\in V_S$.
Then there exists an $\omega\in V_S$ such that
\begin{enumerate}
\item[\rm (i)] $\alpha q(\omega)=1$;
\item[\rm (ii)] the polynomial $\Phi_{\omega}(t)$
(defined in \ref{phidefrel}) is separable of degree
$2n-1$.
\end{enumerate}
\end{lem}

\begin{proof}
According to our notation
we write $V_F$ for $V/\m V$ and
$V_E$ for $V_S/\m V_S$.
The element $\overline{\alpha} \in E$ is a primitive
element of $E/F$. It satisfies the relation
$\overline{\alpha} \overline{q}(\overline{v})=\bar{1}$ for
the vector $\overline{v}\in V_E$.
By Lemma \ref{uslem} applied to $\overline{\alpha}$
and $\overline{v}$ there exists
$\omega'\in V_E$ satisfying the conditions
$\overline{\alpha}\overline{q}(\omega')=\bar{1}$ and
$\spr{\overline{v},\omega'} -\bar{1} \in E^*$,
where $\spr{\, ,\,} : V_E\times V_E \to E$
is the bilinear form associated with
$\overline{\alpha}\overline{q}$.

Now apply Lemma \ref{imlem} to the quadratic space
$(V_S, \alpha q_S)$ and the vectors $v\in V_S$,
$\omega'\in V_E$. We find a vector $\omega \in V_S$
such that
\begin{enumerate}
\item[\rm (i)] $\alpha q(\omega)=1$ in $S$;
\item[\rm (ii)] $\overline{\omega}=\omega'$ in $V_E$.
\end{enumerate}
Property (ii) implies that the reduction
modulo $\m$
of the polynomial $\Phi_{\omega}(t)\in R[t]$ coincides
with the polynomial $\Phi_{\omega'}(t) \in F[t]$.
The polynomial $\Phi_{\omega'}(t)$ is separable
of degree $2n-1$ by  property (iii) of \ref{uslem}.
Hence, $\Phi_{\omega}(t)$ is separable of degree
$2n-1$ (see \ref{seppol}).
\end{proof}

\begin{lem}\label{easylemrel}
Let $(V,q)$ be a quadratic space over $R$.
Then the group of squares
$(R^*)^2$ is contained in $D_q(R)$.
In particular, $a\in D_q(R)$
if and only if $a^{-1}\in D_q(R)$.
\end{lem}

\begin{proof} For any $b\in S^*$ and $q(u)\in R^*$
we have $b^2=q(ub)/q(u)\in D_q(R)$.
\end{proof}

\begin{prop}\label{proprel}
Let $\alpha\in S$ be a primitive invertible
element of the
finite \'etale extension $S/R$.
In particular,
the ring $S$ can be written as $S=R[t]/(f_\alpha(t))$,
where $f_\alpha(t)$ is the minimal polynomial
of $\alpha$.
Assume $\alpha q(v)=1$ for a vector $v \in V_S$.
Then $N^S_R(q(v))\in D_q(R)$.
\end{prop}

\begin{proof}
By Lemma \ref{keylem} there exists 
$\omega\in V_S$ such that
\begin{enumerate}
\item[\rm (i)] $\alpha q(\omega)=1$;
\item[\rm (ii)] the polynomial $\Phi_{\omega}(t)$ 
is separable of degree $2n-1$.
\end{enumerate}
We have $\Phi_{\omega}(\alpha)=\alpha q(\omega)-1=0$.
This implies
$\Phi_{\omega}(t)=c\cdot h(t)f_\alpha(t)$,
where $c\in R^*$ and
$h(t)$ is some monic polynomial of degree $n-1$.
The polynomial $h(t)$ is separable, 
since $\Phi_{\omega}(t)$ is separable.
Clearly, we have the relation 
\begin{equation}\label{polrel}
1 + c\cdot h(t)f_\alpha(t)=tq(\omega(t)).
\end{equation}

The proof proceeds by the induction
on degree of the extension $S/R$.
The case $n=1$ is obvious.
Assume that the proposition holds for
all finite \'etale extensions of degree
strictly less than $n$.

Consider
the finite \'etale extension $T=R[t]/(h(t))$ over $R$,
where  $h(t)$ is the
polynomial appearing in (\ref{polrel}).
Observe that the degree of $T/R$ is $n-1$.
Let $\beta$ be the image of $t$ under the quotient map
$R[t]\to R[t]/(h(t))$. Observe that $\beta$ is
a primitive element of the $R$-algebra $T$ with
the  minimal polynomial $h(t)$.

Consider the reduction modulo the ideal
$(h(t))$ of the relation (\ref{polrel}).
We get $1=\beta q(u)$ in $T$ for some $u\in V_T$.
By Lemma \ref{easylemrel} we get $\beta \in D_q(T)$.
Substituting $t=0$ in (\ref{polrel})  we get
$f_\alpha(0)=-1/(c\cdot h(0))$.
Together with the fact that
for the extension $S=R[t]/(f_\alpha(t))$ over $R$,
$N^S_R(\alpha)=(-1)^{\deg(S/R)}f_\alpha(0)$,
this implies the following chain of relations in $R$:
$$
N^S_R(q(v))=1/N^S_R(\alpha)=(-1)^n /f_\alpha(0)=
c \cdot (-1)^{n-1}h(0)
=c \cdot N^T_R(\beta).
$$
Since $c$ is the leading coefficient of
the polynomial $tq(\omega(t))$ of degree $2n-1$,
it is represented by the form $q$.
Namely, $c=q(\omega_{n-1,0},\ldots,\omega_{n-1,m-1})$,
where $\omega_{n-1,j}\in R$
are the leading coefficients of polynomials $\omega_j(t)$.
By the induction hypothesis the norm
$N^T_R(\beta)$ lies in $D_q(R)$.
Hence $N^S_R(q(v))\in D_q(R)$.
This completes the proof of Proposition \ref{proprel}
\end{proof}

\begin{lem}\label{densqrel}
Consider the multiplicative group $\gm_{,S}$
over the semi-local scheme $\spec S$.
Let $W$ be an open subset of $\gm_{,S}$ such that
for each closed point $x\in \spec S$
the fiber $W_x$ over $x$ is non-empty.
Then there exists
$b\in S^*$ such that $b^2\in W(S)$.
\end{lem}

\begin{proof}
For each closed point $x$ of $\spec S$ there is
an element $a_x$ in the residue field of $x$ such that
$a_x^2\in W_x$.
Since $S$ is semi-local there is an element
$b\in S^*$ such that
$b_x=a_x$ in the residue field of $x$, for each $x$.
Clearly, $b^2\in W(S)$.
\end{proof}

\begin{proof}[Proof of Theorem \ref{normprrel}]
Let $a=q(u) \in S^*$ for certain  $u\in V_S$.
Consider the non-empty open subset $P$
of $R_{S/R}(\gm_{,S})$ defined in \ref{notetrel}.
Clearly, each closed fiber
of $aP$ over $\spec S$ is non-empty.
By Lemma \ref{densqrel}
there exists an element $b\in S^*$ such that
$b^2\in a P(S)$.
It means that $\alpha=a^{-1}b^2$
is in $P(S)$, i.e., primitive.
Replacing $a$ by $\alpha$ and $u$ by
$v=u\cdot b^{-1}$, we get $\alpha q(v)=1$.
  Then, by Lemma \ref{easylemrel} and
Proposition \ref{proprel}, we get the desired
inclusion
$N^S_R(a)=N^S_R(q(v))\cdot N^S_R(b)^2 \in D_q(R)$.
\end{proof}

\begin{ntt}
As before let
$R$ be a semi-local Noetherian domain with infinite residue fields
of characteristic different from 2.
Let $S/R$ be a finite \'etale $R$-algebra of degree $n$.
Let $(V,q)$ be a quadratic space over $R$ of rank $m$.
By $D^0_q(S)$ (resp. $D^1_q(S)$) we denote the set of all even (resp. 
odd) products
of invertible elements of $S$ represented by $q_S$, i.e.,
$$
D^i_q(S)=\{\prod_{j=0}^lq(v_j) \mid v_j\in S^m,\,q(v_j)\in S^*,\, 
l\equiv i\, {\rm mod}\, 2 \},\quad i=0,1.
$$
Observe that $D^0_q(S)$ is a subgroup of the group $D_q(S)$ and 
$(S^*)^2\subset D^0_q(S)$.
Clearly, if $c\in D^i_q(S)$, $i=0,1$, and $b\in D^0_q(S)$,
then $cb\in D^i_q(S)$.
The following result is an obvious consequence of the proof of 
Theorem \ref{normprrel}.
\end{ntt}

\begin{thm}\label{normpreven}
Let $N^S_R:S^*\to R^*$ be the norm map. Then
$$
N^S_R(D^0_q(S)) \subset D^0_q(R),
$$

\end{thm}

\begin{proof}
For a positive integer $n$ we set $D^n_q(S)=D^0_q(S)$ if $n$ is even
and $D^n_q(S)=D^1_q(S)$ if $n$ is odd.
By the proof of Proposition \ref{proprel}
and Theorem \ref{normprrel} it follows that $N^S_R(a)\in D^n_q(R)$,
where $a$ is an element represented by $q_S$ and $S/R$ is an 
extension of degree $n$.
\end{proof}

%
%

\section{Grothendieck's conjecture for  the spinor group}\label{grothconj}

Let $R$ be a local domain with residue field
of characteristic different from 2
and $q$ be a  quadratic space over $R$.
Following \cite[IV.6]{Kn} we define the spinor group (scheme) $Spin_q$
to be $Spin_q(R)=\{x\in S\Gamma_q(R)\mid x\sigma(x)=1\}$,
where $\sigma$ is the canonical involution,
$S\Gamma_q(R)=\{c\in C_0(V,q)^* \mid cVc^{-1}\subset V\}$ is
the special Clifford group and
$C_0(V,q)$ is the even part of the Clifford algebra of
the respective quadratic space $(V,q)$ over $R$.
The present section is devoted to the proof of the following result:
\begin{thm}\label{groth}
Let $R$ be a local regular ring containing an infinite field
of characteristic different from 2. Let $K$ be its quotient field.
Let $q$ be a  quadratic space over $R$.
Then the induced map on the sets of principal homogeneous spaces
$$
H^1_{et}(R,Spin_q)\to H^1_{et}(K,Spin_q)
$$
has trivial kernel,
where $Spin_q$ is the spinor group for the quadratic space $q$.
\end{thm}

\begin{ntt}
Observe that the theorem is a particular case
of Grothendieck's conjecture on principal homogeneous spaces
\cite{Gr},
which states that, for a smooth reductive group scheme $G$ over $R$,
the induced map $H^1_{et}(R,G)\to H^1_{et}(K,G)$ has trivial kernel.
\end{ntt}

\begin{proof}
The proof is based on the results of  \cite{Oj},
\cite{Za1} and \cite{Za2}.

Assume $R$ is a local regular ring containing
a field of characteristic different from 2.
Let $K$ be its quotient field.
We have the following commutative diagram (see \cite[IV.8.2.7]{Kn}):
$$
\xymatrix{
SO_q(R) \ar[r]^{SN} \ar[d] & R^*/(R^*)^2 \ar[r] \ar[d] &
H^1_{et}(R, Spin_q) \ar[r] \ar[d] &  H^1_{et}(R,SO_q) \ar[d] \\
SO_q(K) \ar[r]^{SN} & K^*/(K^*)^2 \ar[r] &
H^1_{et}(K, Spin_q) \ar[r] &  H^1_{et}(K,SO_q),
}
$$
where $SN: SO_q(R) \to H^1_{et}(R,\mu_2)=R^*/(R^*)^2$ is the spinor norm.
The main result of  \cite{Oj} says that the vertical arrow on the right
hand side has trivial kernel (see also \cite[3.4]{Za1}).
Thus, in order to show that the middle one
has trivial kernel, it is enough to check that the induced map on the
cokernels $coker(SN)(R)\to coker(SN)(K)$ is injective.
First, we prove the geometric case.
To do that we use the following slight modification
of the main result of section 2 of \cite{Za1}.

\begin{prop}\label{inject}
Let $R$ be a local regular ring of geometric type over an infinite
field.
Let $F$ be a presheaf from the category of affine schemes over $R$
to abelian groups. Assume $F$ satisfies the axioms
{\rm C}, {\rm E} of \cite{Za1} and the weak versions of the axioms
{\rm TE}, {\rm TA}, {\rm TB} of \cite{Za1}
(see \ref{WTE}). Then the canonical map $F(R)\to F(K)$ is injective.
\end{prop}

\begin{ntt} \label{WTE}
The weak versions of the axioms TE, TA and TB state
\begin{itemize}
\item
For a finite \'etale $R$-algebra $T$ (instead of finitely generated
projective considered in \cite{Za1}) there is given a transfer map
$Tr_{R}^T: F(T) \to F(R)$.
\item
This transfer map is additive in the sense of TA of \cite{Za1}.
\item
For a finitely generated projective $R[t]$-algebra $S$ such that
the algebras $S/(t)$ and $S/(t-1)$ are finite \'etale over $R$
there is the commutative diagram of the axiom TB of \cite{Za1}
$$
\xymatrix{
F(S)\ar[r]\ar[d] & F(S/(t))\ar[d]^-{Tr} \\
F(S/(t-1))\ar[r]^-{Tr} & F(R)
}
$$
\end{itemize}
The proof of \ref{inject} follows immediately after one replaces
the Geometric Presentation Lemma \cite[10.1]{OP} used in section 1.1
of \cite{Za1} by its stronger (\'etale) version
\end{ntt}

\begin{lem}[\'Etale Geometric Presentation Lemma]
Let $R$ be a local essentially smooth algebra over an infinite field $k$,
$m$ its maximal ideal and $S$ an essentially smooth $k$-algebra
which is an integral domain and finite over the polynomial algebra $R[t']$.
Suppose that $e: S\to R$ is an $R$-augmentation and let $I=\ker e$.
Assume that $S/mS$ is smooth over  the residue field $R/m$
at the maximal ideal $e^{-1}(m)/mS$.
Then, given a regular function $f\in S$ such that $S/(f)$ is finite over $R$,
we can find a $t\in I$ such that
\begin{itemize}
\item
$S$ is finite over $R[t]$;
\item
There is an ideal $J$ comaximal with $I$ and such that $I\cap J=(t)$;
\item
$(f)$ and $J$ are comaximal; \ \   $(f)$ and $(t-1)$ are comaximal;
\item
$S/(t)$ is \'etale over $R$; \ \   $S/(t-1)$ is \'etale over $R$.
\end{itemize}
\end{lem}

\begin{proof} See \cite[6.1]{Za2}
\end{proof}

Consider the presheaf of abelian groups
$F:T\mapsto coker(SN)(T)$.
According to \ref{inject}
to prove the mentioned injectivity
we have to show that the functor $F$ satisfies axioms C, E and
the weak versions of axioms TE, TA, TB  of \cite{Za1} (see \ref{WTE}).
Axioms C, E, TA and TB hold by the same arguments as in
sections 3.2 and 3.4 of \cite{Za1}.

Consider, for instance, the proof of axiom E. First, following the proof of E.(a) and E.(b)
of \cite[3.2]{Za1}
for a given quadratic space $q$ over $R$
we construct the $R$-algebra $\tilde{S}$, two
quadratic spaces $q_1=q\otimes_A \tilde{S}$ and
$q_2=q\otimes_R \tilde{S}$ over $\tilde{S}$
and the isomorphism $\Psi$ between them such that
the restrictions $q_1|_R$ and $q_2|_R$
coincide with $q$ and the restriction $\Psi|_R$
is the identity.
Then following the proof of E.(c)
the isomorphism $\Psi$ induces a commutative diagram
$$
\xymatrix{
SO_{q_1}(T) \ar[r]^{\Psi}\ar[d]^{SN} &
SO_{q_2}(T)\ar[d]^{SN} \\
T^*/(T^*)^2 \ar[r]^{\Psi}& T^*/(T^*)^2
}
$$
for any $\tilde{S}$-algebra $T$.
Taking the cokernels of the vertical arrows
we obtain the desired functor
transformation $\Phi: F_1(T)\to F_2(T)$ of E.(c).

Hence, in order to prove the injectivity,
it remains to produce
a well-defined transfer map $Tr_R^S: F(S) \to F(R)$
for any finite \'etale extension $S/R$ of a local regular ring
of geometric type over an infinite field.

To produce such a map it suffices to take the norm map $N^S_R:S^*\to R^*$
and to check the inclusion
$$
N^S_R(SN(SO_q(S)))\subset SN(SO_q(R)).
$$
Since in the semi-local case $SN(SO_q(S))=D^0_q(S)$ and $SN(SO_q(R))=D^0_q(R)$
(see \cite[IV.6]{Kn},
\cite[III.3.21]{Ba}) it remains to check the inclusion
$$
N^S_R(D^0_q(S)) \subset D_q^0(R).
$$
This inclusion holds by Theorem \ref{normpreven}.
Hence, the norm map $N^S_R:S^* \to R^*$ induces the desired transfer map
$Tr^S_R: F(S) \to F(R)$. This completes the proof of Theorem \ref{groth}
in the geometric case.

Finally, to extend our result
to the case of a local regular
ring $R$ containing an infinite field of characteristic different from 2
we use Popescu's approximation theorem
\cite[7.5]{Za2}. We refer to the item 1 of section 5 of \cite{Za1}
for the precise arguments.
\end{proof}

\newpage

\noindent
Manuel Ojanguren \\
Ecole Polytechnique F\'{e}d\'{e}rale de Lausanne \\
Institut de G\'eom\'etrie, Alg\`{e}bre et Topologie \\
B\^atiment BCH \\
1015 Lausanne \\
Switzerland \\
e-mail: manuel.ojanguren@epfl.ch

\

\noindent
Ivan Panin \\
St. Petersburg Department of \\
V.A. Steklov Institute of Mathematics \\
27, Fontanka \\
St. Petersburg 191023 \\
Russia \\
e-mail: panin@pdmi.ras.ru

\

\noindent
Kirill Zainoulline \\
Fakult\"at f\"ur Mathematik \\
Universit\"at Bielefeld \\
Postfach 100131 \\
33501 Bielefeld \\
Germany \\
e-mail: kirill@mathematik.uni-bielefeld.de

\end{document}